\documentclass[mathpazo]{cicp}
\usepackage{booktabs}
\usepackage{url}

\begin{document}
\title{Strang-type preconditioners for solving fractional diffusion equations by boundary value methods}


 \author[X.-M. Gu et.~al.]{Xian-Ming Gu\affil{1}\comma\corrauth,
       Ting-Zhu Huang\affil{1}, Xi-Le Zhao\affil{1}, Hou-Biao Li\affil{1}, Liang Li\affil{1}}
 \address{\affilnum{1}\ School of Mathematical Sciences,
          University of Electronic Science and Technology of China,
          Sichuan 611731, P. R. China.
          }
 \emails{{\tt guxianming@yahoo.cn} (X.-M. Gu), {\tt tingzhuhuang@126.com} (T.-Z. Huang)}

\begin{abstract}
The finite difference scheme with the shifted Gr\"{u}nwarld formula is employed to semi-discrete the fractional diffusion equations. 
This spatial discretization can reduce to the large system of ordinary differential equations (ODEs) with initial values. Recently, 
boundary value method (BVM) was developed as a popular algorithm for solving large systems of ODEs. This method requires the solutions of one or more 
nonsymmetric, large and sparse linear systems. In this paper, the GMRES method with the block circulant preconditioner is proposed for solving these 
linear systems. One of the main results is that if an $A_{\nu_1,\nu_2}$-stable boundary value method is used for an $m$-by-$m$ system of ODEs, 
then the preconditioner is invertible and the preconditioned matrix can be decomposed as $I+L$, where $I$ is the identity matrix and the rank of $L$ is 
at most $2m(\nu_1+\nu_2)$. It means that when the GMRES method is applied to solve the preconditioned linear systems, the method will converge in at most $2m(\nu_1+\nu_2)+1$ iterations. 
Finally, extensive numerical experiments are reported to illustrate the effectiveness of our methods for solving the fractional diffusion equations.
\end{abstract}


\maketitle

\section{Introduction}
\label{sec1}
\quad During recent years, the concept of fractional derivatives, and their applications to modelling anomalous diffusion phenomena are widely recognised by engineers and mathematicians. Fractional diffusion equations are useful for applications in which a cloud of particles spreads faster than predicted by the classical equation. FDEs arise in research topics
including modeling chaotic dynamics of classical conservative systems \cite{GMZ}, turbulent flow \cite{BAC,MFS}, groundwater
contaminant transport \cite{DBSW1,DBSW2}, and applications in
biology \cite{RLM}, finance \cite{EOL,MRE}, image processing \cite{JBX,JMB}, hydrology \cite{BBD}
and other physics issues \cite{IMS}. For example, anomalous diffusion is a possible mechanism underlying plasma
transport in magnetically confined plasmas, and the fractional order space derivative operators can
be used to model such transport mechanism. As there are very few cases of FDEs in which the closed-form
analytical solutions are available, numerical solutions for FDEs become main ways and then have been developed intensively, such as (compact) finite difference method \cite{TAML,FLVA,MMM1,MMM2,MMM3,DAMI}, finite element method \cite{WHDF,VJE,KBN,JPRC}, discontinuous Galerkin method \cite{XJHZT,KMWM,WHDJS} and other numerical methods \cite{CPE,KXWH,CMCFL,AAZCO,LJS,BBMK}.

However, due to the nonlocal character of the fractional differential operator, it was shown that a
naive discretization of the FDE, even though implicit, leads to unconditionally unstable
\cite{MMM2,MMM3}. Moreover, most numerical methods for FDEs tend to generate full coefficient matrices, which require computational cost of $\mathcal{O}(N^3)$ and storage of $\mathcal{O}(N^2)$, where $N$ is the
number of grid points \cite{HWKX2}. It is quite different from second-order diffusion equations which
usually yield sparse coefficient matrices with $\mathcal{O}(N)$ nonzero entries and can be solved very
efficiently by fast iterative methods with $\mathcal{O}(N)$ complexity.

To overcome the difficulty of the stability, Meerschaet and Tadjeran \cite{MMM2,MMM3} proposed a
shifted Gr\"{u}nwald discretization to approximate FDEs. Their method has been proven to be
unconditionally stable. Later, Wang, et. al \cite{HWKX2} discovered that the full coefficient
matrix by the Meerschaet-Tadjeran¡¯s method holds a Toeplitz-like structure. More precisely,
such a full matrix can be written as the sum of diagonal-multiply-Toeplitz matrices. Thus
the storage requirement is significantly reduced from $\mathcal{O}(N^2)$ to $\mathcal{O}(N)$. It is well known
that the matrix-vector multiplication for the Toeplitz matrix can be computed by the fast
Fourier transform (FFT) with $\mathcal{O}(N \log N)$ operations \cite{RCMN1,RCMN2,MNI}. With this advantage, Wang
and Wang \cite{KXWH} employed the conjugate gradient normal residual (CGNR) method to solve
the discretized system of the FDE by the Meerschaet-Tadjeran's method. Thanks to the
Toeplitz-like structure, the cost per iteration by the CGNR method is of $\mathcal{O}(N \log N)$. The
convergence of the CGNR method is fast with smaller diffusion coefficients \cite{KXWH} (in that case
the discretized system is well-conditioned). Nevertheless, if the diffusion coefficient functions
are not small, the resulting system will become ill-conditioned and hence the CGNR method
converges very slowly. To overtake this shortcoming, Pang and Sun \cite{HPHS} proposed a multigrid
method to solve the discretized system of the FDE by the Meerschaet-Tadjeran's method.
With the damped-Jacobi method as the smoother, the multigrid algorithm can preserve the
computational cost per iteration as $O(N \log N)$ operations. Numerical results showed that
their multigrid method converges very fast, even for the ill-conditioned systems. However,
from the theoretical point of view, the linear convergence of their multigrid method, despite
a very simple case (both diffusion coefficients are equal and constant), has not been proven,
see \cite{HPHS} for details. Recently, Lei and Sun \cite{SLLH} proposed a robust CGNR method with the circulant preconditioner to solve FDEs by the Meerschaet-Tadjeran's method under the conditions that the diffusion coefficients are constant and the ratio is bounded away from zero. The convergence analysis of  their method can be archived more easily than the multigrid method does.

In this paper, we firstly induce the FDEs to be a system of ODEs by the spatial discretization method (semi-discretization). In particular, we apply the GMRES \cite{YSMS} with the block-circulant type preconditioners for solving linear systems arising from the application of BVMs, which is a relatively new method based on the linear multistep formulae to solve ODEs. Boundary value methods (BVMs) are unconditionally stable and are high-accuracy schemes for solving initial value problems (IVPs) based on the linear multistep formulas \cite{LBDT,PGPM,PAFMDT}. Unlike Runge-Kutta or other initial value methods (IVM), BVMs achieve the advantage of both good stability and high-order accuracy \cite{PAFMDT,LBDTS2}. The main purpose of this paper is to investigate the effectiveness of preconditioning technique on the speed of the resulting iterative processes of boundary value methods for solving FDEs.

The paper is organized as follows. In Section 2, the background of the spital discretization
for the FDE to reduce the system of ODEs is reviewed. Then we introduce that how to result in the linear systems by block-BVMs. In Section 3, we construct the block circulant-type preconditioner and BCCB preconditioner. Then the invertibility of two different kinds of preconditioner and the convergence rate and operation cost of the preconditioned GMRES method are also studied.
In Section 4, extensive numerical results are reported to demonstrate the efficiency of the proposed method.


\section{Semi-discretization for FDEs and boundary value methods}
\quad In this paper, we study an initial-boundary value problem of the FDE as follows,
\begin{equation}
\label{eq1.1}
\left \{ {\begin{array}{l}
\frac{\partial u(x,t)}{\partial t}=d_{+}(x,t)\frac{\partial^{\alpha}u(x,t)}{\partial_{+}x^{\alpha}}+
d_{-}(x,t)\frac{\partial^{\alpha}u(x,t)}{\partial_{-}x^{\alpha}}+f(x,t),\\
x\in(x_L, x_R),\qquad t\in(t_0,T],\\
u(x_L,t)=u(x_R,t)=0,\qquad 0\leq t\leq T,\\
u(x,t_0)=u_0(x),\qquad x\in [x_L,x_R],\\
\end{array}}
\right.
\end{equation}
where $\alpha\in(1, 2)$ is the order of the fractional derivative, $f(x, t)$ is the source (or sinks) term, and
diffusion coefficient functions $d_{\pm}(x,t)$ are nonnegative; i.e., $d_{\pm}(x,t)\geq0,\ d_{+}(x,t)+d_{-}(x,t)\neq0$. The function $u(x,t)$ can be interpreted as representing the concentration of a particle plume undergoing anomalous diffusion.
\subsection{FDM semi-discretization for FDEs}
\label{sec2}
\quad Meerschaert and Tadjeran \cite{MMM3} have shown that using the shifted Gr\"{u}nwald formula to approximate the two-sided fractional derivatives of order $\alpha\in (1,2)$ leads to stable numerical schemes.
We begin this method, it is known that the left-sided and the right-sided fractional derivatives $\frac{\partial^{\alpha}u(x,t)}{\partial_{+}x^{\alpha}}$ and
$\frac{\partial^{\alpha}u(x,t)}{\partial_{-}x^{\alpha}}$ are defined in the Gr\"{u}nwald-Letnikov form \cite{IPF}
\[
\frac{\partial^{\alpha}u(x,t)}{\partial_{+}x^{\alpha}}=\lim_{\Delta x\rightarrow 0^{+}}
\frac{1}{\Delta x^{\alpha}}\sum_{k=0}^{\lfloor(x-x_{L})/\Delta x \rfloor}g_{k}^{(\alpha)}u(x-k\Delta x,t),
\]
\[
\frac{\partial^{\alpha}u(x,t)}{\partial_{-}x^{\alpha}}=\lim_{\Delta x\rightarrow 0^{+}}
\frac{1}{\Delta x^{\alpha}}\sum_{k=0}^{\lfloor(x_{R}-x)/\Delta x \rfloor}g_{k}^{(\alpha)}u(x+k\Delta x,t),
\]
where $\lfloor\cdot\rfloor$ denotes the floor function, and $g_{k}^{(\alpha)}$ is the alternating fractional binomial coefficient
given as
\begin{equation}
\label{eq1.2}
\left \{ {\begin{array}{l}
g_{0}^{(\alpha)}=1,\\
g_{k}^{(\alpha)}=\frac{(-1)^k}{k!}\alpha(\alpha-1)\cdots(\alpha-k+1),\quad\quad k=1,2,3,\ldots,
\end{array}}
\right.
\end{equation}
which can be evaluated by the recurrence relation
\[
g_{k+1}^{(\alpha)}=(1-\frac{\alpha+1}{k+1})g_{k}^{(\alpha)},\quad\quad k=0,1,2,\ldots,
\]
Let $N$ be positive integers and $\Delta x = \frac{x_R-x_L}{N+1}$ be the sizes of spatial
grid. We define a spatial and temporal partition $x_i = x_L + i\Delta x$ for $i = 0, 1,\ldots,N + 1$. Let $u_i = u(x_i, t), d_{\pm,i} =
d_{\pm}(x_i, t)$, and $f_i = f(x_i, t)$. The shifted Gr\"{u}nwald approximation in \cite{MMM2,MMM3} is as
follows,
\[
\frac{\partial^{\alpha}u(x_i,t)}{\partial_{+}x^{\alpha}}=
\frac{1}{\Delta x^{\alpha}}\sum_{k=0}^{i+1}g_{k}^{(\alpha)}u_{i-k+1},
\]
\[
\frac{\partial^{\alpha}u(x_i,t)}{\partial_{-}x^{\alpha}}=
\frac{1}{\Delta x^{\alpha}}\sum_{k=0}^{N-i+2}g_{k}^{(\alpha)}u_{i+k-1},
\]
where $g_{k}^{(\alpha)}$ is defined in (\ref{eq1.2}), and the corresponding spital semi-discretization for the (\ref{eq1.1}), i.e., its corresponding systems of ODEs as follows,
\begin{equation}\label{eq1.3}
\left \{ {\begin{array}{l}
\frac{d\mathbf{u}(t)}{dt}=J_N\mathbf{u}(t)+\mathbf{f}(t),\qquad t\in (t_0,T],\\
\mathbf{u}(t_0)=[u_0(x_1),u_0(x_2),\ldots,u_0(x_N)]^{\mathrm{T}}=\mathbf{u}_0,
\end{array}}
\right.
\end{equation}
where $\mathbf{u}(t) = [u_1, u_2,\ldots, u_N]^{\mathrm{T}},\ \mathbf{f}(t) = [f_1, f_2,\ldots, f_N]^{\mathrm{T}},\ \Delta x^{\alpha}=\frac{(x_R-x_L)^{\alpha}}{(N+1)^{\alpha}}$, and $J_N$ be the coefficient matrix with an appropriate size and can be written in the following
\begin{equation}\label{eq1.4}
J_N=\frac{1}{\Delta x^{\alpha}}[D_{+}G_{\alpha}+D_{-}G_{\alpha}^{\mathrm{T}}],
\end{equation}
with $D_{\pm}=\mathrm{diag}(d_{\pm,1},d_{\pm,2},\ldots,d_{\pm,N})$ and
\begin{equation}
\label{eq1.5}
G_{\alpha}=
\left[\begin{array}{cccccc}
g_{1}^{(\alpha)}&g_{0}^{(\alpha)}&0&\cdots&0&0\\
g_{2}^{(\alpha)}&g_{1}^{(\alpha)}&g_{0}^{(\alpha)}&0&\cdots&0\\
\vdots&g_{2}^{(\alpha)}&g_{1}^{(\alpha)}&\ddots&\ddots&\vdots\\
\vdots&\ddots&\ddots&\ddots&\ddots&0                          \\
g_{N-1}^{(\alpha)}&\ddots&\ddots&\ddots&g_{1}^{(\alpha)}&g_{0}^{(\alpha)}\\
g_{N}^{(\alpha)}&g_{N-1}^{(\alpha)}&\cdots&\cdots&g_{2}^{(\alpha)}&g_{1}^{(\alpha)}
\end{array}
\right]_{N\times N}
\end{equation}
It is obvious that $G_{\alpha}$
is a Toeplitz matrix (see [8, 9]). Therefore, it can be stored with $N +1$
entries \cite{HWKX2}. Furthermore, the matrix-vector multiplication for the Toeplitz-like matrix $J_N$
in (6) can be obtained in $\mathcal{O}(N \log N)$ operations by the FFT; see \cite{RCMN1,HPHS}. The alternating fractional binomial coefficient $g_{k}^{(\alpha)}$ have some useful properties, that are observed in \cite{MMM2,MMM3,HWKX2}, and are summarized in the following proposition.
\begin{proposition}(\cite{SLLH})
Let $1 <\alpha < 2$ and $g^{(\alpha)}_k$ be defined in (\ref{eq1.2}). We have
\begin{equation}\label{eq1.6}
\left \{ {\begin{array}{l}
g_{0}^{(\alpha)}=1,\ \ g_{1}^{(\alpha)}=-\alpha<0,\ \ g_{2}^{(\alpha)}>g_{3}^{(\alpha)}>\cdots>0, \\
\sum\limits_{k=0}^{\infty}g_{k}^{(\alpha)}=0,\ \ \sum\limits_{k=0}^{\infty}g_{k}^{(\alpha)}<0,\ \forall n\geq 1.
\end{array}}
\right.
\end{equation}
\end{proposition}
Also, we give the following conclusion which is very essential for theoretical analysis,
\begin{proposition}
Let $1 <\alpha < 2$ and $g^{(\alpha)}_k$ be defined in (\ref{eq1.2}). All eigenvalues of $J_N$ fall inside the open disc
\[
\{ z\in \mathbb{C}:\ |z-\gamma_i|<-\gamma_i\},\quad i=1,\ldots,N
\]
where $\gamma_i=(r_{+,i}+r_{-,i})g_{1}^{(\alpha)}<0$ are constants.
\end{proposition}
\textbf{Proof}. Here the entries of matrix $J_N$ are given by
\begin{equation}
\label{eq1.10t}
p_{ij}=
\left
\{\begin{array}{l}
(r_{+,i}+r_{-,i})g_{1}^{(\alpha)},\qquad j=i,\\
r_{+,i}g_{2}^{(\alpha)}+r_{-,i}g_{0}^{(\alpha)}, \quad j=i-1, \\
r_{+,i}g_{0}^{(\alpha)}+r_{-,i}g_{2}^{(\alpha)}, \quad j=i+1,\\
r_{+,i}g_{i-j+1}^{(\alpha)},\qquad \qquad j<i-1,\\
r_{-,i}g_{j-i+1}^{(\alpha)},\qquad \qquad j>i+1,
\end{array}
\right.
\end{equation}
where $r_{\pm,i}=\frac{d_{\pm,i}}{\Delta x^{\alpha}}\geq 0$. It is not hard to find that $p_{ij}\leq0$ for all $i\neq j$, then all the Gershgorin disc of the matrix are centered at $\gamma_i=(r_{+,i}+r_{-,i})g_{1}^{(\alpha)}<0$ with radius
\[
\begin{array}{l}
R_i=\sum\limits_{j=1,j\neq i}^{N}|p_{ij}|=r_{+,i}\sum\limits_{k=0,k\neq 1}^{i}g_{k}^{(\alpha)}+r_{-,i}\sum\limits_{k=0,k\neq 1}^{N-i+1}g_{k}^{(\alpha)}
\ \ (r_{+,i}+r_{-,i}\neq 0)\\
\qquad\qquad\qquad\ \ < (r_{+,i}+r_{-,i})\sum\limits_{k=0,k\neq 1}^{\infty}g_{k}^{(\alpha)}=-(r_{+,i}+r_{-,i})g_{1}^{(\alpha)}=-\gamma_i,
\end{array}
\]
by the properties of the sequence $g_{k}^{(\alpha)}$; see Proposition 2.1.
\begin{remark}
It is worth to note that:
\begin{enumerate}
\item[(\romannumeral1)] The real parts of all eigenvalue of the matrix $J_N$ are strictly negative for all $N$.

\item[(\romannumeral2)] The absolute values of all eigenvalues of the matrix $J_N$ are bounded above by $\max\limits_{1\leq j\leq N}\{2|\gamma_i|\}$.
\end{enumerate}
\end{remark}
\subsection{Boundary value methods (BVMs)}
\quad Next, we discuss a class of robust numerical methods called the boundary value methods (BVMs) for solving the systems of ODEs, see \cite{LBDT,AOHX}. Using the $\mu$-step block-BVM over By using a $v$-step LMF over a uniform mesh
\[t=t_0+jh,\quad j=0,1,\ldots,s,\]
where $h = (T -t_0)/s$ is the step size for the discretization of Eq. (\ref{eq1.3}), we obtain
\begin{equation}\label{eq1.7}
\sum_{i=-\nu}^{\mu-\nu}\alpha_{i+\nu}\mathbf{u}^{(k+1)}_{n+i}=h\sum_{i=-\nu}^{\mu-\nu}\beta_{i+\nu}\mathbf{g}_{n+i},\quad n=\nu,\ldots,s-\mu+\nu
\end{equation}
Here, $\mathbf{u}_n$ is the discrete approximation to $\mathbf{u}(t_n)$, $\mathbf{g}_n=J_N\mathbf{u}_n+\mathbf{f}_n$, and $\mathbf{f}_n=\mathbf{f}(t_n)$. Also, Eq. (\ref{eq1.7}) requires $\nu$ initial conditions and $\mu-\nu$ final conditions which are provided by the following $\mu-1$ additional equations:
\begin{equation}\label{eq1.8}
\sum_{i=0}^{\mu}\alpha_{i}^{(j)}\mathbf{u}_{i}=h\sum_{i=0}^{\mu}\beta_{i}^{(j)}\mathbf{g}_{i},\quad j=1,2,\ldots,\nu-1,
\end{equation}
and
\begin{equation}\label{eq1.9}
\sum_{i=0}^{\mu}\alpha_{\mu-i}^{(j)}\mathbf{u}_{n-i}=h\sum_{i=0}^{\mu}\beta_{\mu-i}^{(j)}\mathbf{g}_{n-i},\quad j=s-\mu+\nu+1,\ldots,s.
\end{equation}

The coefficients $\{\alpha_{k}^{(j)}\}$ and $\{\beta_{k}^{(j)}\}$ in Eq. (\ref{eq1.8}) and Eq. (\ref{eq1.9}) should be chosen such that truncation errors in theses $\mu-1$ equations are of the same order as that in Eq. (\ref{eq1.7}).  We combine Eqs. (\ref{eq1.7})-(\ref{eq1.9}) and the initial condition $\mathbf{u}(t_0)=\mathbf{u}_0$, a discrete linear system of Eq. (\ref{eq1.6}) is given by the following block matrix form,
\begin{equation}\label{eq1.10}
M\mathbf{u}\equiv(A\otimes I_N-hB\otimes J_N)\mathbf{u}=\mathbf{e}_1\otimes \mathbf{u}_0 +h(B\otimes I_N)\mathbf{f}.
\end{equation}
where
\[\mathbf{e}_1=(1,0,\ldots,0)^T \in \mathbb{R}^{s+1},\qquad \mathbf{u}=(\mathbf{u}^{T}_0,\ldots,\mathbf{u}^{T}_s)^{T}\in\mathbb{R}^{(s+1)N},\quad
\mathbf{f}=(\mathbf{f}^{T}_0,\ldots,\mathbf{f}^{T}_s)^{T}\in\mathbb{R}^{(s+1)N}.\] In (\ref{eq1.10}), the matrix $A\in \mathbb{R}^{(s+1)\times(s+1)}$ is given by:
\[
A=\left[
\begin{array}{ccccccc}
1&\cdots&0\\
\alpha_{0}^{(1)}&\cdots&\alpha_{\mu}^{(1)}\\
\vdots&\vdots&\vdots\\
\alpha_{0}^{(\nu-1)}&\cdots&\alpha_{\mu}^{(\nu-1)}\\
\alpha_{0}&\cdots&\alpha_{\mu}\\
&\alpha_{0}&\cdots&\alpha_{\mu}\\
&&\ddots&\ddots&\ddots\\
&&&\ddots&\ddots&\ddots\\
&&&&\alpha_{0}&\cdots&\alpha_{\mu}\\
&&&&\alpha_{0}^{(s-\mu+\nu+1)}&\cdots&\alpha_{\mu}^{(s-\mu+\nu+1)}\\
&&&&\vdots&\vdots&\vdots\\
&&&&\alpha_{0}^{(s)}&\cdots&\alpha_{\mu}^{(s)}
\end{array}
\right],
\]
and $B\in \mathbb{R}^{(s+1)\times(s+1)}$ is defined similarly by using $\beta$'s instead of $\alpha$'s in $A$ and the first row of $B$ is zeros.

Usually the resulting linear system (\ref{eq1.10}) is large and ill-conditioned, and solving it is a core problem in the application of BVMs. If a direct solver is applied to solve the system (\ref{eq1.10}), the operation cost can be very high for practical application. Therefore interest has been turned to iterative solvers, such as GMRES method. As we know that a clustered spectrum often translates in rapid convergence of GMRES \cite{MBGHG}, so we use the GMRES method for solving the resulting linear system (\ref{eq1.10}). In order to accelerate the convergence of iterations, we construct some block circulant-type preconditioners.
\section{Construction of preconditioners and convergence analysis}
\quad In this section, we will show how to construct the block circulant-type preconditioners for accelerating the iterative solver and show that these preconditioners are invertible if an $A_{\nu_1,\nu_2}$-stable BVM is used. Meanwhile, some theoretical analyses on both the convergence rate of iterative solver and operation cost are also investigated.
\subsection{Construction of preconditioners}
\quad To mimic the terminology of \cite{RHCXQJYH} and neglecting the perturbations in the upper left and low right corners of $A$ and $B$, we give the first preconditioner for Eq. (\ref{eq1.10}):
\begin{equation}
\label{eq1.11}
S=s(A)\otimes I_{N}-hs(B)\otimes J_N,
\end{equation}
where
\[s(A)=\left[
\begin{array}{ccccccccc}
\alpha_{\nu}&\cdots&\alpha_{\mu}&&&&\alpha_{0}&\ldots&\alpha_{\nu-1}\\
\vdots&\ddots&&\ddots&&&&\ddots&\vdots\\
\alpha_{0}&&\ddots&&\ddots&&&&\alpha_{0}\\
&\ddots&&\ddots&&\ddots&&0\\
&&\ddots&&\ddots&&\ddots\\
&0&&\ddots&&\ddots&&\ddots\\
\alpha_{\mu}&&&&\ddots&&\ddots&&\alpha_{\mu}\\
\vdots&\ddots&&&&\ddots&&\ddots&\vdots\\
\alpha_{\nu+1}&\cdots\alpha_{\mu}&&&&&\alpha_{0}&\cdots&\alpha_{\nu}
\end{array}
\right]
.\]
and $s(B)$ is defined similarly by using $\{\beta_i\}_{i=0}^{\mu}$ instead of $\{\alpha_i\}_{i=0}^{\mu}$ in $s(A)$. The $\{\alpha_i\}_{i=0}^{\mu}$ and $\{\beta_i\}_{i=0}^{\mu}$ here are the coefficients in Eq. (\ref{eq1.7}). We note that $s(A)$ and $s(B)$ are the generalized Strang-type circulant preconditioners of $A$ and $B$ respectively, see \cite{RCMN2}.

Moreover, we also can propose the Strang-type BCCB preconditioner, which can be constructed for solving Eq. (\ref{eq1.10})
\begin{equation}\label{eqx1}
S^{(2)}=(s(A)\otimes I_N-hs(B)\otimes s(J_N))
\end{equation}
for $J_N$ being a Toeplitz-like matrix with structure as the sum of diagonal-multiply-Toeplitz matrices. Here we define the $s(J_N)$ as following form
\[
s(J_N)=\frac{1}{\Delta x^{\alpha}}(\overline{d}_{+}s(G_{\alpha})+\overline{d}_{-}s(G^{\mathrm{T}}_{\alpha}))
\]
with $\overline{d}_{\pm}=\frac{1}{N}\sum\limits_{i=1}^{N}d_{\pm,i}$. More precisely, the first columns of $s(G_{\alpha})$ and $s(G^{\mathrm{T}}_{\alpha})$ are given by
\[\left[
\begin{array}{c}
g_{1}^{(\alpha)}\\
\vdots\\
g_{\lfloor\frac{N+1}{2}\rfloor}^{(\alpha)}\\
0\\
\vdots\\
0\\
g_{0}^{(\alpha)}
\end{array}
\right]\qquad \mathrm{and}\qquad
\left[
\begin{array}{c}
g_{1}^{(\alpha)}\\
g_{0}^{(\alpha)}\\
0\\
\vdots\\
0\\
g_{\lfloor\frac{N+1}{2}\rfloor}^{(\alpha)}\\
g_{0}^{(\alpha)}
\end{array}
\right]
.\]
As we know, Lei and Sun \cite{SLLH} proposed Strang circulant preconditioner to approximate the coefficient matrix matrix with structure as the sum of diagonal-multiply-Toeplitz matrices. The convergent behavior of this method are very efficient and robust in numerical experiments. So we take the similar strategy to construct the preconditioner $S^{(2)}$ (\ref{eqx1}). The advantage of BCCB preconditioners is that the operation cost in each iteration of Krylov subspace methods for the preconditioned systems is much less than that required by using any block-circulant preconditioners.

Next, we will display that the preconditioner $S$ is invertible provided that the given BVM is stable and the eigenvalue of $J_m$ are in the negative half of the complex plane $\mathbb{C}$. Also the invertibility of the preconditioner $S^{(2)}$ will be analyzed and improved. The stability of a BVM is related to two characteristic polynomials of degree $\mu$, defined as follows:

\begin{equation}\label{eq1.12}
\rho(z)=z^{\nu}\sum_{j=-\nu}^{\mu-\nu})\alpha_{j+\nu}z^{j}\qquad \mathrm{and}\qquad\sigma(z)=z^{\nu}\sum_{j=-\nu}^{\mu-\nu})\beta_{j+\nu}z^{j}.
\end{equation}

\begin{definition}(\cite[p. 101]{LBDT})
Consider a BVM with the characteristic polynomials $\rho(z)$ and $\sigma(z)$ given by (\ref{eq1.11}). The region
\[
\begin{array}{l}
\mathcal{D}_{\nu,\mu-\nu}=\{q\in \mathbb{C}:\rho(z)-q\sigma(z)\ has\ \nu\ zeros\ inside\ |z|=1\\
\qquad\qquad\qquad\ \ \ \ and\ \mu-\nu\ zeros\ outside\ |z|=1\}
\end{array}
\]
is called the region of $A_{\nu,\mu-\nu}$-stability of the given BVM. Moreover, the BVM is said to be $A_{\nu,\mu-\nu}$-stable if
\[
\mathbb{C}^{-}\equiv\{q\in\mathbb{C}:\ \mathrm{Re}(q)<0\}\subseteq \mathcal{D}_{\nu,\mu-\nu}.
\]
\end{definition}
\begin{theorem}(\cite{RHCXQJYH})
If the BVM for (\ref{eq1.3}) is $A_{\nu,\mu-\nu}$-stable and $h\lambda_k(J_N)\in \mathcal{D}_{\nu,\mu-\nu}$ where $\lambda_k(J_N)\ (k=1,\ldots,N)$ are the eigenvalues of $J_N$, then the preconditioner $S$ in (\ref{eq1.11}) is invertible.
\end{theorem}

In particular, we have
\begin{corollary}(\cite{RHCXQJYH})
If the BVM for (\ref{eq1.3}) is $A_{\nu,\mu-\nu}$-stable and $h\lambda_k(J_N)\in \mathbb{C}^{-}$, then the preconditioner $S$ is invertible.
\end{corollary}

In fact, we can find that the eigenvalues of $J_N$ are in the negative half of the complex plane $\mathbb{C}^{-}$ by the Proposition 2.2 and Remark 2.1. So if we add the condition that the given BVM is stable, we can immediately conclude that the preconditioner $S$ is invertible. It means that this preconditioner can be expected to be robust and efficient.

Similar to Theorem 3.1, we can show that if the BVM for (1) is $A_{\nu,\mu-\nu}$-stable and the eigenvalues of $s(J_N)$ satisfy
\[\lambda_k(s(J_N))\in \mathbb{C}^{-}\]
for $k = 1,\ldots,N$, then the preconditioner $S^{(2)}$ is invertible.

However, for some special FDEs problem, the matrix $J_N$ is usually full Toeplitz-like structure, but $s(J_N)$ may be singular. Note that the eigenvalues of $S^{(2)}$ are given by
\begin{equation}\label{eq1.13}
\lambda_{jk}(S^{(2)}) = \phi_j -h\psi_j\lambda_k(s(J_N)),\quad j = 0,\ldots,s,\quad k = 1,\ldots,N,
\end{equation}
where $\phi_j$ and $\varphi_j$ are eigenvalues of $s(A)$ and $s(B)$ respectively. When some eigenvalues of $s(J_N)$ are zero,
then some eigenvalues of $S^{(2)}$ is the same as the eigenvalues of the matrix $s(A)$. It is well-known that the eigenvalues of
the circulant matrix $s(A)$ can be expressed as the following sum, see [11],
\[\phi_j =\sum_{r=-\nu}^{\mu-\nu} \alpha_{r+\nu}\omega^{rj},\quad \omega= e^{2\pi \mathbf{i}/(s+1)},\quad j = 0,\ldots,s,
\]
where $\alpha_{r+\nu}$ are given by (6).

From the characteristic polynomials defined in (11), the coefficients must satisfy the
consistent conditions,
\[
\rho(1) = 0\quad \mathrm{and}\quad \rho'(1) = \sigma(1).
\]
Thus, we have
\[\phi_0 =\rho(1) = 0\]
for any consistent BVM. From (15), we know that $S^{(2)}$ is singular when some eigenvalues of $s(J_N)$ are zero. In this case, we move the zero eigenvalue of $s(A)$ to a nonzero value.
More precisely, we change the matrix $s(A) = F\mathrm{diag}(\phi_0,\ldots,\phi_s) F^{*}$ to
\[\widetilde{s}(A)\equiv F\mathrm{diag}(\widetilde{\phi}_0,\ldots,\phi_s)F^{*},\]
where $\widetilde{\phi}_0\equiv \mathrm{Re}(\phi_s)$ and $F$ is the Fourier matrix. Define
\begin{equation}\label{eqxc}
\widetilde{S}^{(2)}\equiv \widetilde{s}(A)\otimes I_N - hs(B) \otimes s(J_N),
\end{equation}
we can also prove that $\widetilde{S}^{(2)}$ is invertible, see \cite{SLLXQ} for a detail.

From the conclusions of \cite{SLLH}, we can obtain the following theorem,
\begin{theorem}
All eigenvalues of circulant matrices $s(G_{\alpha})$ and $s(G^{\mathrm{T}}_{\alpha})$ fall inside the open disc
\[
\{ z\in \mathbb{C}:\ |z+\alpha|<\alpha\}.
\]
\end{theorem}
\textbf{Proof.} The proof of this theorem is greatly similar to that of the Lemma 1 in \cite{SLLH}, we omit here.

By the theorem, we can find that the parts of all eigenvalues of $s(G_{\alpha})$ and $s(G^{\mathrm{T}}_{\alpha})$ are strictly negative for all $N$. Moreover, we know that $\overline{d}_{\pm}\geq0,\ \overline{d}_{+}+\overline{d}_{-}\neq 0$. So we can conclude that
\[
\mathrm{Re}(\lambda_k(s(J_N)))=\overline{d}_{+}\mathrm{Re}(s(G_{\alpha}))+\overline{d}_{-}\mathrm{Re}(s(G^{\mathrm{T}}_{\alpha}))<0.
\]
It means that the eigenvalues of $s(J_N)$ are in the negative half of the complex plane $\mathbb{C}^{-}$ and then both the preconditioners $S^{(2)}$ and $\widetilde{S}^{(2)}$ are invertible provided that the given BVM is stable by the Theorem 2 of \cite[p. 32]{SLLXQ}.
\subsection{Convergence rate and operation cost}
\quad As the statements in \cite{RHCXQJYH}, we have the following theorems for the convergence rates,
\begin{theorem}(\cite{RHCXQJYH})
We have
\[
S^{-1}M = I + L
\]
where $I$ is the identity matrix and the rank of $L$ is at most $2m\mu$. Therefore, when the GMRES method is applied to solving $S^{-1}M\mathbf{y} = S^{-1}\mathbf{b}$, the method will converge in at most $2m\mu + 1$ iterations in exact arithmetic.
\end{theorem}

Lei and Jin \cite{SLLXQ} proved that when $J_N$ is a Toeplitz matrix in the \textit{Wiener class} \cite{RCMN1,RCMN2}, the preconditioned matrix $(\widetilde{S}^{(2)})^{-1}M$ can be written as the sum of the identity matrix, a matrix with rank $\mathcal{O}(N)$, a matrix with rank $\mathcal{O}(s)$ and a matrix with small norm. In fact, for Eq. (\ref{eq1.3}), when we take
\begin{equation}\label{eqfg}
d_{+,i}=d_{+}\geq 0,\quad  d_{-,i}=d_{-}\geq 0 \quad \mathrm{and}\ \ d_{+}+d_{-}\neq 0,
\end{equation}
for all $i=1,\ldots,N$. Then we obtain a nonsymmetric Toeplitz matrix as follows,
\begin{equation}\label{eqgu}
\begin{array}{l}
J_N=T_N=\frac{1}{\Delta x^{\alpha}}(d_{+}G_{\alpha}+d_{-}G^{\mathrm{T}}_{\alpha})\\
\quad \ =\frac{1}{\Delta x^{\alpha}}\left[
\begin{array}{ccccc}
d_{+}g^{(\alpha)}_{1}+d_{-}g^{(\alpha)}_{1}&d_{+}g^{(\alpha)}_{0}+d_{-}g^{(\alpha)}_{2}&d_{-}g^{(\alpha)}_{3}&\cdots&d_{-}g^{(\alpha)}_{2}\\
d_{+}g^{(\alpha)}_{2}+d_{-}g^{(\alpha)}_{0}&d_{+}g^{(\alpha)}_{1}+d_{-}g^{(\alpha)}_{1}&\ddots&\ddots&\vdots\\
d_{+}g^{(\alpha)}_{3}&\ddots&\ddots&\ddots&\vdots\\
\vdots&\ddots&\ddots&\ddots&d_{+}g^{(\alpha)}_{0}+d_{-}g^{(\alpha)}_{2}\\
d_{+}g^{(\alpha)}_{N}&\cdots&\cdots&d_{+}g^{(\alpha)}_{2}+d_{-}g^{(\alpha)}_{0}&d_{+}g^{(\alpha)}_{1}+d_{-}g^{(\alpha)}_{1}
\end{array}
\right]\\
\quad \ =[t_{j-k}]_{N\times N}.
\end{array}
\end{equation}
we introduce the \textit{generating function} of the sequence of Toeplitz matrices $\{T_N\}_{N=1}^{\infty}$ \cite{RCMN1}:
\begin{equation}
p(\theta)= \sum_{k=-\infty}^{\infty}t_{k}e^{\mathbf{i}k\theta},
\end{equation}
where $t_k$ is the $k$-th diagonal of $T_N$. The generating function $p(\theta)$ is in the Wiener class if and only if
\[
\sum_{k=-\infty}^{\infty}|t_k|<\infty.
\]
For $T_N$ defined in (\ref{eqgu}), we have
\begin{equation}
p(\theta)= \sum_{k=-\infty}^{\infty}t_{k}e^{\mathbf{i}k\theta}=\sum_{k=-1}^{\infty}g^{(\alpha)}_{k+1}(d_{+}e^{\mathbf{i}k\theta}+d_{-}e^{\mathbf{i}k\theta}).
\end{equation}
and obtain the following theorem
\begin{theorem}
Let $p$ be the generating function of $\{T_N\}_{N=1}^{\infty}$, we conclude that $p$ is in the Wiener class.
\end{theorem}
\textbf{Proof}. By the properties of the sequence $\{g^{(\alpha)}_{k+1}\}_{k=0}^{\infty}$ given in (\ref{eq1.6}), we have
\[
\begin{array}{l}
\sum\limits_{k=-\infty}^{\infty}|t_k|=\frac{1}{\Delta x^{\alpha}}(d_{+}+d_{-})\sum\limits_{k=-1}^{\infty}|g^{(\alpha)}_{k+1}|\\
\quad \quad \quad\ \ \ =\frac{1}{\Delta x^{\alpha}}(d_{+}+d_{-})\Big(-2g^{(\alpha)}_{1}+\sum\limits_{k=0}^{\infty}g^{(\alpha)}_{k}\Big)\\
\quad \quad \quad\ \ \ =\frac{2\alpha}{\Delta x^{\alpha}}(d_{+}+d_{-})<\infty.
\end{array}
\]
Thus $p$ is in the Wiener class.

Therefore, when we take the assumption of (\ref{eqfg}), then we can say that the preconditioned matrix $(\widetilde{S}^{(2)})^{-1}M$ can be written as the sum of the identity matrix, a matrix with rank $\mathcal{O}(N)$, a matrix with rank $\mathcal{O}(s)$ and a matrix with small norm by the use of Theorem 3 in \cite[p. 34]{SLLXQ}. As a consequence, the spectrum of $(\widetilde{S}^{(2)})^{-1}M$ is clustered around $(1,0)\in \mathbb{C}$. Moreover, the GMRES method, when applied for solving the preconditioned linear system
\[
(\widetilde{S}^{(2)})^{-1}M\mathbf{y}=(\widetilde{S}^{(2)})^{-1}\mathbf{b}
\]
will converge fast. Therefore, a detailed analysis for the convergence rate could be carried out in the future work.

Regarding the cost per iteration, the main work in each iteration for the GMRES method is the matrix-vector multiplication
\[
S^{-1}M\mathbf{z}= (s(A)\otimes I_N- hs(B)\otimes J_N)^{-1}(A\otimes I_N-hB\otimes J_n)\mathbf{z}
\]
where $\mathbf{z}$ is a vector, see for instant Saad \cite{YSMS}. Since $A$ and $B$ are band matrices and $J_N$ and $ s(B)$ is a full matrix, the matrix-vector multiplication $M\mathbf{z}=(A\otimes I_N-hB\otimes J_n)\mathbf{z}$ can be implemented not slowly.

To calculate $S^{-1}M\mathbf{z}$, since $s(A)$ and $s(B)$ are circulant matrices, we have the following decompositions by fast Fourier transform (FFT)
\[
s(A)= F\Lambda_{A}F^{*}\quad\ \mathrm{and}\quad\ s(B)= F\Lambda_{B}F^{*}
\]
where $\Lambda_A$ and $\Lambda_B$ are diagonal matrices containing the eigenvalues of $s(A)$ and $s(B)$ respectively, see [5]. It follows that
\[
S^{-1}(M\mathbf{z}) = (F^{*}\otimes I_N)(\Lambda_{A} \otimes I_N -h\Lambda_{B} \otimes J_N)^{-1}(F\otimes I_N)(M\mathbf{v}).
\]
This product can be obtained by using FFT and solving $s$ (Toeplitz-like) linear systems of order $m$. It follows that the total number of operations per iteration is $\pi_1 ms\log s+\pi_2smn$, where $n$ is the number of nonzeros of $J_N$, and $\pi_1$ and $\pi_2$ are some positive constants. For
comparing the computational cost of the method with direct solvers for the linear system (\ref{eq1.10}), we refer to [10]. However, in the case of numerical method for FDEs, the coefficient matrix $J_N$ is full, it means that $n$ is much large. We need to take much time to solve $s$ (Toeplitz-like) linear systems of order $m$, this shortage will keep the preconditioner $S$ from becoming the efficient one. In order to overcome this shortage, we propose the preconditioners $S^{(2)}$ and $\widetilde{S}^{(2)}$. For simplicity, we assume that $s+1=N$ in the following analysis of the operation cost of preconditioners $S^{(2)}$ and $\widetilde{S}^{(2)}$. Regarding the cost in each iteration of the GMRES method, the main work is the matrix-vector multiplication
\[
(\widetilde{S}^{(2)})^{-1}M\mathbf{v}\equiv \widetilde{s}(A)\otimes I_N-hs(B) \otimes s(J_N))^{-1}M\mathbf{v},
\]
where $\mathbf{v}$ is a vector. Since $(\widetilde{S}^{(2)})^{-1}$ can be diagonalized by the 2-dimensional Fourier matrix, i.e.,
\[
(\widetilde{S}^{(2)})^{-1}M\mathbf{v}\equiv(F_{s+1}\otimes F_N)(\Lambda_{A}\otimes I_N- h\Lambda_{B}\otimes \Lambda_{J_N})^{-1}(F^{*}_{s+1} \otimes F^{*}_N)(M\mathbf{v}).
\]
where $s(J_N)=F^{*}_{N}\Lambda_{J_N}F_{N}$ and $\Lambda_{J_N}$ is a diagonal matrix holding the eigenvalues of $s(J_N)$. The matrix-vector multiplication $(\widetilde{S}^{(2)})^{-1}M\mathbf{v}$ can be obtain within $\mathcal{O}(N^2\log N)$ operations by using the FFT. For the Strang-type block-circulant preconditioner $S$ defined as the form (\ref{eq1.11}), in each iteration, there are m Toeplitz-like systems of order $m$ needed to be solved. Thus, the complexity in each iteration of the preconditioners $S^{(2)}$ and $\widetilde{S}^{(2)}$ is much lower.
\section{Numerical experiments}
\quad In this section, we solve two different FDE problems (\ref{eq1.1}) numerically by the BVM and the GMRES
method together with the circulant-type preconditioners in Section 2-3. We also compare the Strang-type BCCB preconditioners $S^{(2)}$ and $\widetilde{S}^{(2)}$ with the Strang-type block-circulant preconditioner $S$. Number of iterations required for convergence and CPU time of those methods are reported. In these examples, the BVM we used here is the fifth order GAM which has $\mu=4$. Its formulae and the additional initial and final conditions can be found in Ref. \cite{LBDT}.

All experiments are performed in {\sf MATLAB} 2011b and all the computations are done on a Inter(R) Pentium(R) CPU 2.80GHz PC with 3.85G available memory. We use the {\sf MATLAB}-provided M-file `gmres'
(see {\sf MATLAB} on-line documentation) to solve the preconditioned systems. We use donations ``ITS" and ``CPU" to represent the number of iterations and CPU elapsed time (mean value from ten times repeated experiments) of implementing GMRES(20) solver, respectively. In our tests, the initial guess is the zero vector and
stopping criterion in the GMRES method is
\[
\frac{\|\mathbf{r}_q\|_2}{\|\mathbf{r}_0\|_2}<10^{-8},
\]
where $\mathbf{r}_q$ is the residual after the $q$th iterations.
%
%

\textbf{Example 4.1}. In this example, we solve the initial-boundary value problem of FDE (1) with
source term $f(x, t)\equiv 0$, for order of fractional derivatives $\alpha = 1.2$ and $1.5$. The spatial
domain is $[x_L, x_R] = [0, 2]$ and the time interval is $[0, T] = [0, 1]$. The initial condition $u(x, 0)$
is the following Gaussian pulse
\[
u(x, 0) = \exp\Big(-\frac{(x-x_c)^2}{2\xi^2}),\qquad x_c=1.2,\quad \xi=0.08,
\]
and the diffusion coefficients
\[
d_{+}(x,t) \equiv0.6,\quad\ \mathrm{and}\quad\ d_{-}(x,t)\equiv 0.5.
\]

\begin{table}[htbp]\footnotesize
\begin{center}
\begin{tabular*}{14.3cm}{p{20pt}p{21pt}p{30pt}p{35pt}p{35pt}p{32pt}p{32pt}p{30pt}p{28pt}p{28pt}}
\multicolumn{10}{p{280pt}}{\small \textbf{Table 1}}\\
\multicolumn{10}{p{280pt}}{\small The number of iterations and CPU time (s) of GMRES(20) solver for Example 1 with $\alpha=1.2$.}\\
\hline &&$\ \ \ \ \ \ \ \ \ \ \ I$&&$\ \ \ \ \ \ \ \ \ \ \ S$&&$\ \ \ \ \  \ \ \ \ S^{(2)}$&&$\ \ \ \ \ \ \ \ \ \
\widetilde{S}^{(2)}$\\
\cline{3-10} $n$ &$s$  &ITS &CPU      &ITS&CPU      &ITS        &CPU              &ITS&CPU \\
\hline
24 &16  &82  &0.0781    &9  &0.0310   &15 &0.0234           &17 &0.0263 \\
   &32  &130 &0.1714    &8  &0.0475   &15 &0.0348           &17 &0.0359  \\
   &64  &265 &0.5934    &8  &0.0935   &15 &0.0521           &17 &0.0588 \\
   &128 &457 &2.1064    &7  &0.2291   &14 &0.0989           &17 &0.1142 \\
48 &16  &174 &0.1562    &9  &0.0783   &19 &0.0308           &20 &0.0336  \\
   &32  &198 &0.3122    &8  &0.1249   &18 &0.0442           &20 &0.0485  \\
   &64  &260 &0.8279    &8  &0.2988   &17 &0.1148           &20 &0.1455 \\
   &128 &460 &2.5592    &7  &0.4213   &17 &0.1363           &20 &0.1948 \\
96 &16  &234 &0.2654    &9  &0.2811   &23 &0.0457           &26 &0.0532 \\
   &32  &262 &0.6425    &8  &0.4654   &23 &0.1172           &26 &0.1314 \\
   &64  &339 &1.3739    &7  &0.7948   &21 &0.1713           &26 &0.1901 \\
   &128 &393 &2.8865    &7  &1.4978   &21 &0.2927           &26 &0.3248 \\
\hline
\end{tabular*}
\end{center}
\end{table} \begin{table}[htbp]\footnotesize                                                                                                        \begin{center}
\begin{tabular*}{14.3cm}{p{20pt}p{21pt}p{30pt}p{35pt}p{35pt}p{32pt}p{32pt}p{30pt}p{28pt}p{28pt}}
\multicolumn{10}{p{280pt}}{\small \textbf{Table 2}}\\
\multicolumn{10}{p{280pt}}{\small The number of iterations and CPU time (s) of GMRES(20) solver for Example 1 with $\alpha=1.5$.}\\
\hline &&$\ \ \ \ \ \ \ \ \ \ \ I$&&$\ \ \ \ \ \ \ \ \ \ \ S$&&$\ \ \ \ \  \ \ \ \ S^{(2)}$&&$\ \ \ \ \ \ \ \ \ \
\widetilde{S}^{(2)}$\\
\cline{3-10} $n$ &$s$  &ITS &CPU      &ITS&CPU      &ITS        &CPU              &ITS&CPU \\
\hline
24 &16  &141  &0.1208    &10 &0.0328   &19 &0.0271           &27 &0.0341  \\
   &32  &201  &0.2607    &9  &0.0558   &19 &0.0407           &27 &0.0522  \\
   &64  &237  &0.5416    &8  &0.0977   &19 &0.0668           &28 &0.0918  \\
   &128 &378  &1.8091    &8  &0.2038   &18 &0.1632           &28 &0.2056  \\
48 &16  &259  &0.2327    &10 &0.0782   &25 &0.0375           &34 &0.0466  \\
   &32  &282  &0.4308    &9  &0.1357   &25 &0.0602           &36 &0.0803  \\
   &64  &313  &0.9682    &8  &0.2703   &25 &0.1498           &37 &0.1902  \\
   &128 &431  &2.3788    &8  &0.4262   &24 &0.2326           &37 &0.3042  \\
96 &16  &394  &0.4289    &10 &0.2802   &23 &0.0574           &47 &0.0807  \\
   &32  &428  &0.9987    &9  &0.4838   &23 &0.1453           &49 &0.1898  \\
   &64  &532  &2.0718    &8  &0.8568   &23 &0.2386           &51 &0.3118  \\
   &128 &632  &4.6136    &8  &1.4922   &32 &0.3788           &51 &0.5459  \\
\hline
\end{tabular*}
\end{center}
\end{table}

Table 1-2 list the number of iterations required for convergence of the GMRES method with different precondtioners and their corresponding CPU time. In the tables, $I$ means no preconditioner is used, and $S,\ S^{(2)}$ and $\widetilde{S}^{(2)}$ denote the Strang-type block-circulant preconditioners, Strang-type and modified Strang-type BCCB preconditioners respectively, see (\ref{eq1.11}), (\ref{eqx1}) and (\ref{eqxc}).

From Example 4.1, the number of iterations of both $S^{(2)}$ and $\widetilde{S}^{(2)}$ are larger than those of Strang-type block-circulant preconditioner $S$. But the operation cost per iteration of both $S^{(2)}$ and $\widetilde{S}^{(2)}$ is less than those of $S$. As we can see from Table 1-2, the CPU time of $S^{(2)}$ is less than those of the others especially when $n$ and $s$ are large. Moreover, the matrix $J_N$ is ill-conditioned when $n$ is large. The performance of $S^{(2)}$ is the best in terms of the CPU time. We strongly suggest that the preconditioner $S^{(2)}$ is a good choice and we do not need to formulate the complete matrix $J_N$ in order to save storage. Especially, when $J_N$ is the Toeplitz-like structure ($d_{\pm}(x,t)\neq const$). In order to further illustrate the effectiveness of the block-circulant preconditioners, we list the spectra of the original matrix $M$ and the preconditioned matrices $S^{-1}M,\ (S^{(2)})^{-1}M (\widetilde{S}^{(2)})^{-1}M$ in Fig. 1.

\begin{figure}[htbp]\label{fig:1}
\centering
\includegraphics[width=2.92in,height=2.45in]{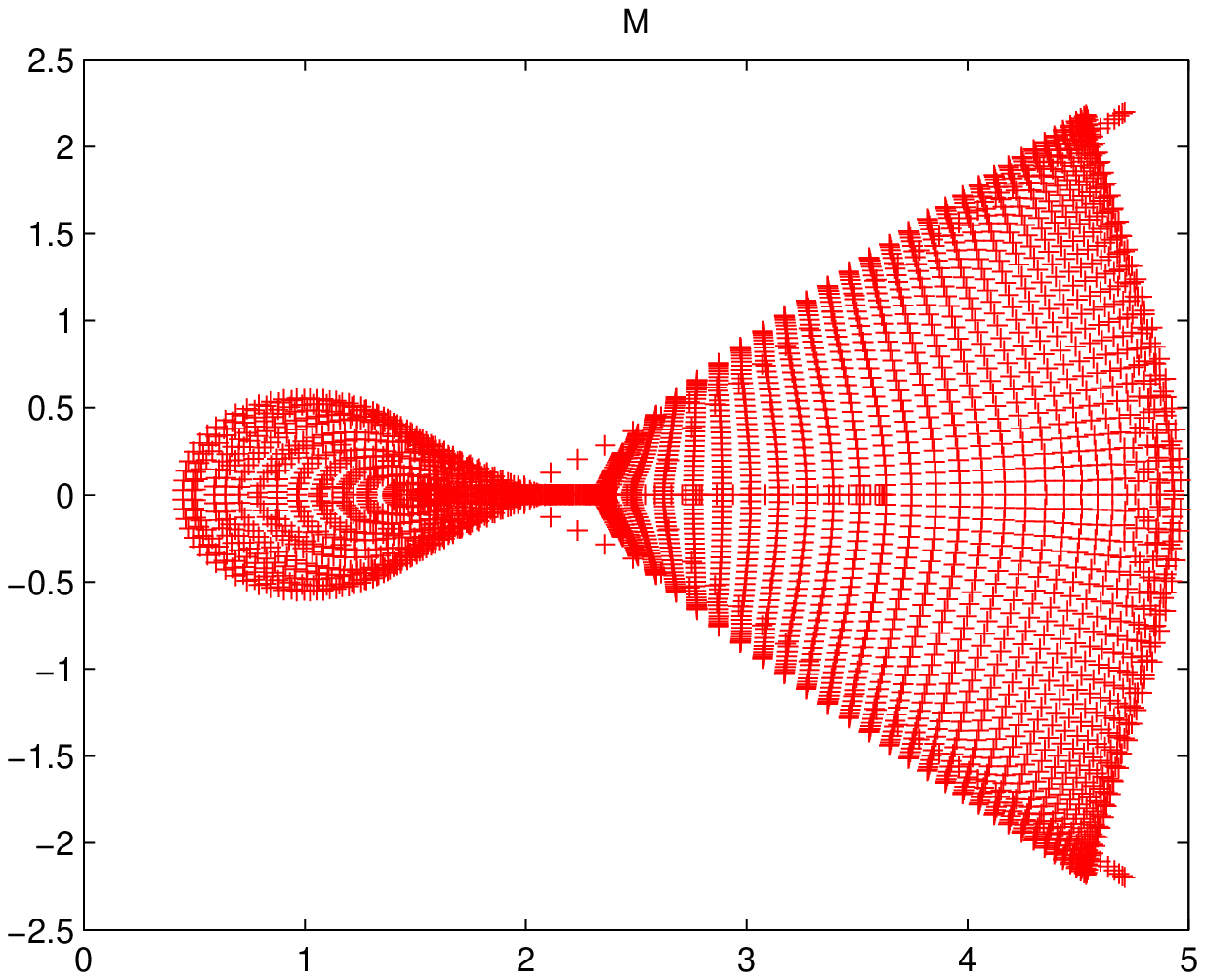}
\includegraphics[width=2.92in,height=2.45in]{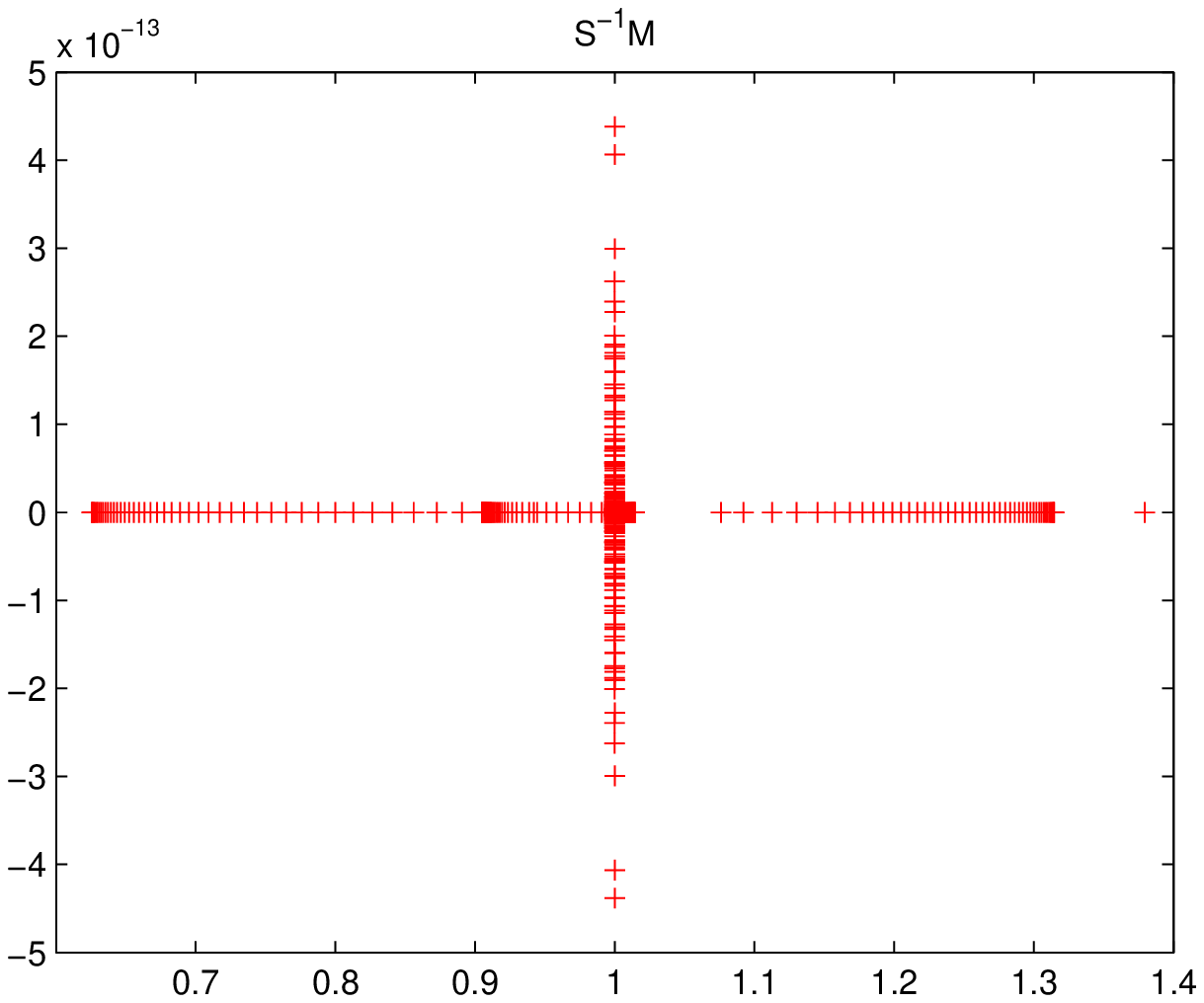}
\includegraphics[width=2.92in,height=2.45in]{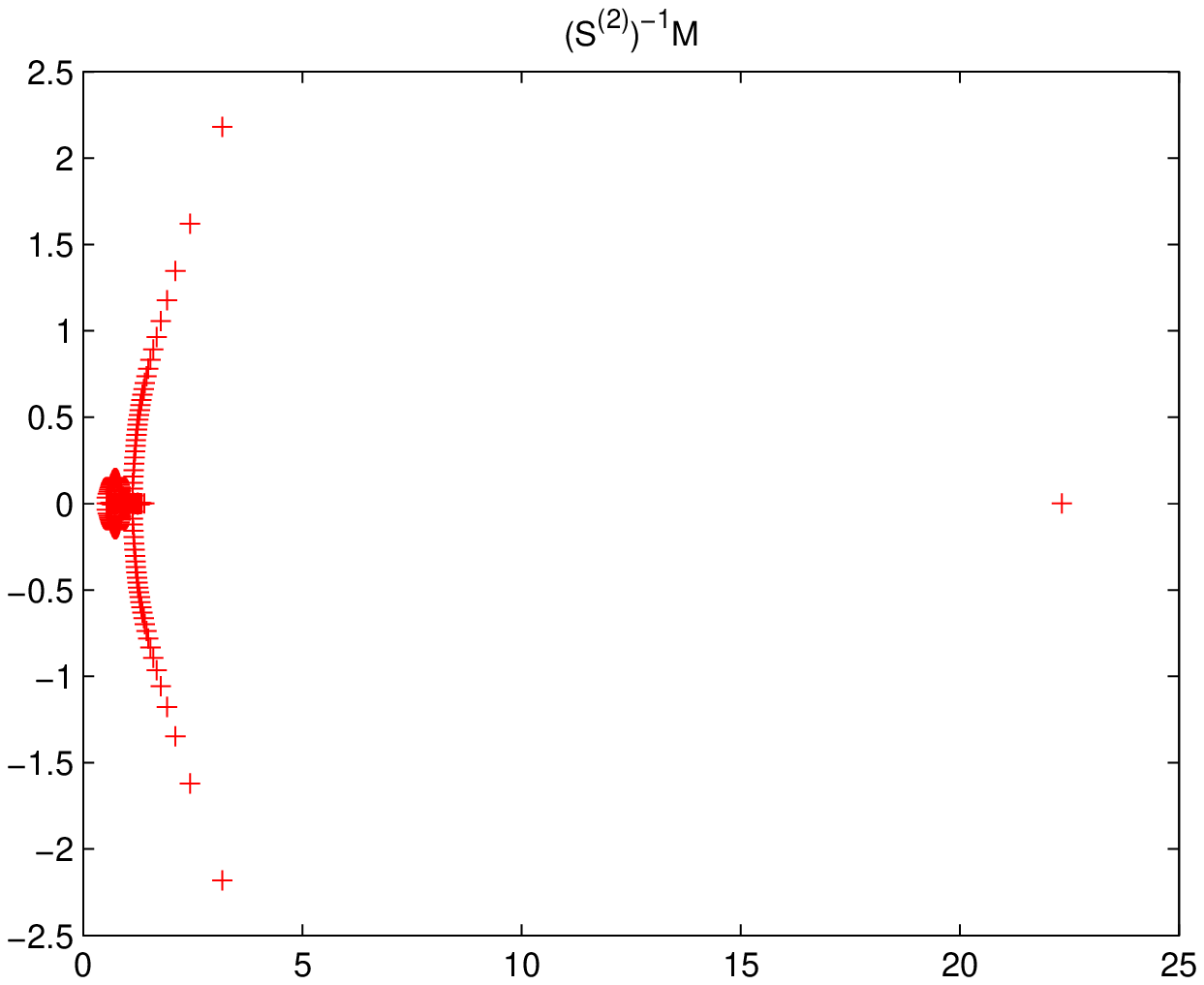}
\includegraphics[width=2.92in,height=2.45in]{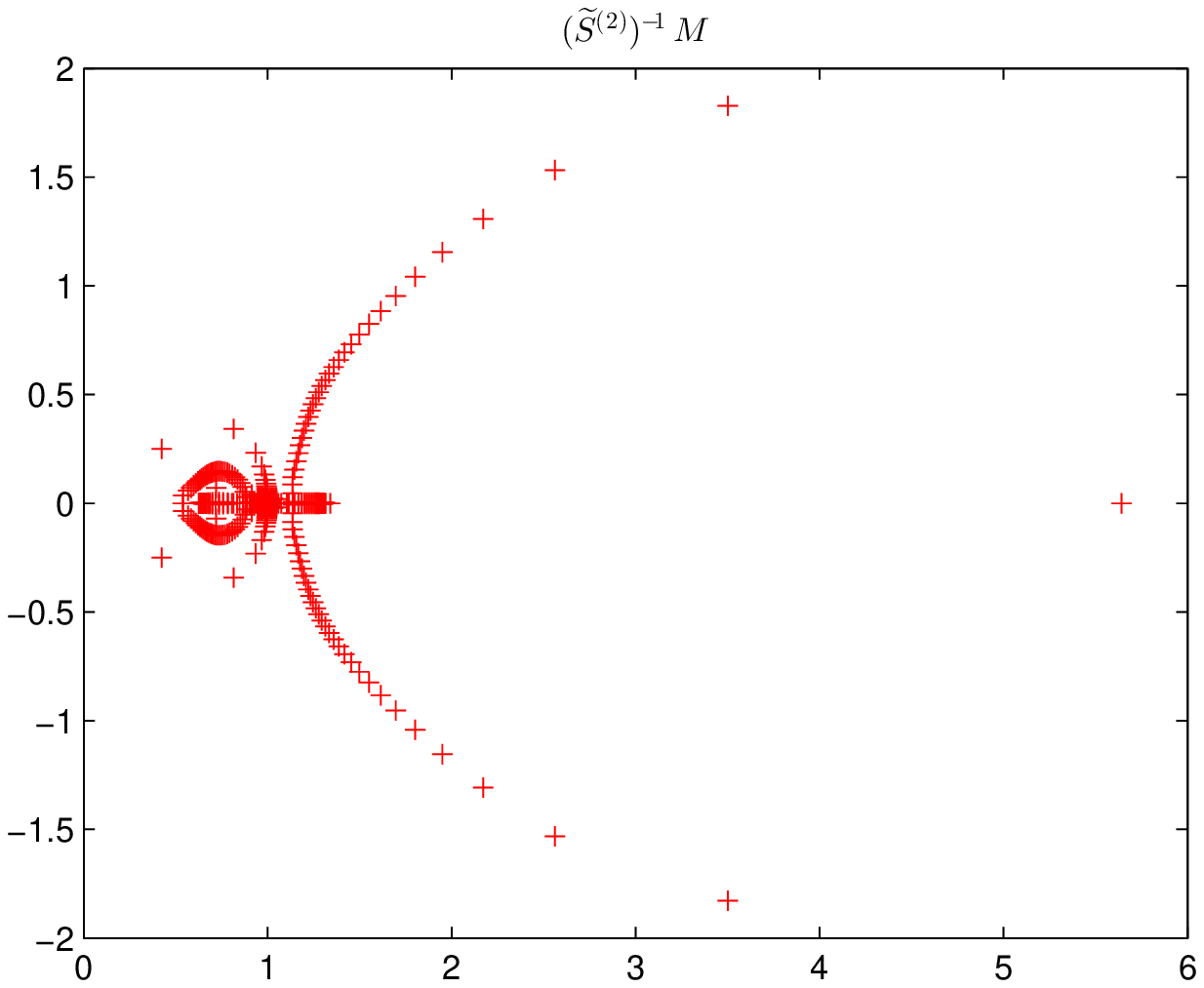}
\caption{\small The spectra of matrix $M$ and different preconditioned matrices with $n=48,\ s=64$ for the Example.}
\end{figure}

\section*{Acknowledgments}
The authors would like to thank Dr. Siu-Lung Lei of the University of Macau for his helpful discussions. This research is supported by NSFC (61170311, 11101071, 1117105, 51175443, 11271001), Chinese Universities Specialized Research Fund for the Doctoral Program (20110185110020)
 and the Fundamental Research Funds for China Scholarship Council.


\end{document}